\documentclass[12pt]{amsart}
\usepackage{amsmath,amsthm,latexsym,amscd,amsbsy,amssymb,amsfonts,fleqn,leqno,
euscript, pb-diagram}
\usepackage{color}

\numberwithin{equation}{section}

\newtheorem{thm}{Theorem}[section]

\newtheorem{conj}[thm]{Conjecture}
\newtheorem{cor}[thm]{Corollary}
\newtheorem{qu}[thm]{Question}

\def\leukfrac#1/#2{\leavevmode
               \kern.1em
                \raise.9ex\hbox{\the\scriptfont0 ${}_#1$}
                \hskip -1pt\kern-.1em
                /\kern-.15em\lower.10ex\hbox{\the\scriptfont0 ${}_#2$}}

\theoremstyle{definition}

\theoremstyle{remark}

\theoremstyle{definition}

\newcommand{\bd}{\mathrm{bd}}

\begin{document}

\title[Homogeneous locally compact spaces]
{Homogeneous locally compact spaces}

\author{V. Valov}
\address{Department of Computer Science and Mathematics,
Nipissing University, 100 College Drive, P.O. Box 5002, North Bay,
ON, P1B 8L7, Canada} \email{veskov@nipissingu.ca}

\date{\today}
\thanks{The author was partially supported by NSERC
Grant 261914-19.}


\begin{abstract}
This is a survey of the recent results and unsolved problems about locally compact homogeneous metric spaces. Mostly, homogeneous finite-dimensional $ANR$-spaces are discussed. 
\end{abstract}
\maketitle\markboth{}{Homogeneous $ANR$,s}





\section{Introduction}

In this paper we survey the most recent results and unsolved problems concerning homogeneous finite-dimensional spaces. It can be considered as an continuation of our previous paper \cite{VV}. As in \cite{VV}, we are mainly interested in locally compact metric homogeneous $ANR$s, but some problems concerning more general homogeneous spaces are also considered. Recall that a space $X$ is {\em homogeneous} if for every two points $x,y\in X$ there is a homeomorphism $h$ mapping $X$ onto itself with $h(x)=y$. This implies that $X$ is {\em locally homogeneous}, i.e. for every two points $x,y\in X$, there 
exists a homeomorphism $h$ mapping a neighborhood $U_x$ of $x$ onto a neighborhood 
$U_y=h(U_x)$ of $y$ and satisfying the condition $h(x_1)=y$.

There are many interesting problems about homogeneous spaces. Probably, the best known is the Bing-Borsuk conjecture \cite{bb} stating that every $n$-dimensional homogeneous metric $ANR$-compactum, $n\geq 3$, is an $n$-manifold. This conjecture is true in dimensions 1 and 2 \cite{bb}. Recently, Bryant-Ferry \cite{BF} provided a revised version of their paper containing counter-examples to that conjecture. They constructed for every $n\geq 6$ infinitely many, topologically distinct, homogeneous $ANR$-compacta of dimension $n$ that are not topological manifolds. So, this conjecture is still open in dimensions 3, 4 and 5. 
Another open problem is whether there is a non-degenerated finite-dimensional locally homogeneous, in particular homogeneous, $AR$-spaces, see \cite{bb}, \cite{bo}.
On the other hand, finite-dimensional locally compact homogeneous $ANR$s share many properties with Euclidean manifolds, see for example \cite{VV1}, \cite{VV}, \cite{vv1}. So, although homogeneous finite-dimensional $ANR$-compacta may not be Euclidean manifolds, it still interesting to what extend they have common properties with Euclidean manifolds. The survey paper of Halverson-Repov\v{s} \cite{hr} contains more information for different type of homogeneity. 

Recall that a metric space $X$ is an absolute neighborhood retract (br., $ANR$) if for every embedding of $X$ as a closed subset of a metric space $M$ there exists a neighborhood $U$ of $X$ in $M$ and a retraction $r:U\to X$, i.e. a continuous map $r$ with $r(x)=x$ for all $x\in X$. Contractible $ANR$-spaces form the class of absolute retracts (br., $AR$).

Unless stated otherwise, all spaces are locally compact separable metric and all maps are continuous. Reduced \v{C}ech homology $\check{H}_n(X;G)$ and cohomology $\check{H}^n(X;G)$ with coefficients from an abelian group $G$ are considered. Singular homology and cohomology groups are denoted, respectively, by $H_k(X;G)$ and $H^k(X;G)$. By a dimension we mean the covering dimension $\dim$, the cohomological dimension with respect to a group $G$ is denoted by $\dim_G$. Recall that $\dim_G$ is the largest integer $n$ such that there exists a closed set $A\subset X$ with $\check{H}^n(X,A;G)\neq 0$.

\section{Homogeneous spaces and the Bing-Borsuk conjecture}

  It is interesting whether some of the counter examples to the Bing-Borsuk conjecture constructed by Bryant-Ferry have the stronger version of homogeneity, the so called strong local homogeneity. Recall that a space $X$ is {\em strongly locally homogeneous} (br., $SLH$) if every point in $X$ has a base of neighborhoods $U$ such that for every $y,z\in U$ there is a homeomorphism $h:X\to X$ with $h(y)=z$ and $h(x)=x$ for all $x\not\in U$.
Every strongly locally homogeneous space $X$ is homogeneous provided $X$ is connected. Moreover, if no two-point set disconnetes a connected strongly locally homogeneous space $X$, then $X$ is {\em $n$-homogeneous} for all $n\geq 1$ \cite{ba}: if $A, B$ are two $n$-elements subsets of $X$, then there is a homeomorphism $h$ on $X$ such that $h(A)=B$. 
The question whether there is an $SLH$-counter example to the Bing-Borsuk conjecture is interesting because every Euclidean manifold has this property.
So, the Bing-Borsuk conjecture can be restated:
\begin{qu}
Is it true that every $n$-dimensional strongly locally homogeneous $ANR$-compactum is an $n$-manifold? 
\end{qu}

Jakobsche \cite{ja} proved that the $3$-dimensional Bing-Borsuk conjecture implies the Poincare conjecture. Assuming the Poincare conjecture is not true, he constructed a 3-dimensional homogeneous compact $ANR$-space which is not a manifold. Any such an example have the additional property of $n$-homogeneity for all $n$, see \cite[Theorem 8.1]{ja1}. There is a strong expectation that Jacobsche's construction provides a 3-dimensional $SLH$ compact $ANR$-space which is not a manifold. Therefore, we have another natural question:
\begin{qu}
Is it true that the restated Bing-Borsuk conjecture in dimension 3 imply the Poincare conjecture? 
\end{qu}

Topological $n$-manifolds $X$ have the following property: For every $x\in X$ the groups $H_k(X,X\setminus\{x\};\mathbb Z)=0$ if $k<n$   and $H_n(X,X\setminus\{x\};\mathbb Z)=\mathbb Z$. A space with this property is said to be a {\em $\mathbb Z$-homology $n$-manifold}. A {\em generalized $n$-manifold} is a locally compact $n$-dimensional $ANR$-space which is a $\mathbb Z$-homology $n$-manifold.
Every generalized $(n\leq 2)$-manifold is known to be a topological $n$-manifold \cite{w}. On the other hand, for every $n\geq 3$ there exists a generalized $n$-manifold $X$ such that $X$ is not locally Euclidean at any point, see for example \cite{can}. Let us mention that the Bryant-Ferry \cite{BF} counter examples to the Bing-Borsuk conjecture are generalized n-manifolds with the disjoint disks property, $n\geq 6$. 

Bryant \cite{br1} suggested another modification of the Bing-Borsuk conjecture:
\begin{conj}\cite{br1}$[$Modified Bing-Borsuk conjecture$]$
Every locally compact homogeneous $ANR$-space of dimension $n\geq 3$ is a generalized $n$-manifold.
\end{conj}
A partial result concerning the Modified Bing-Borsuk conjecture is an old result of Bredon \cite{bre}, reproved by Bryant \cite{br2}:
\begin{thm}\cite{bre},\cite{br2}
If $X$ is a locally compact homogeneous $ANR$-space of dimension $n$ such that the groups $H_k(X,X\setminus\{x\};\mathbb Z)$, $k\leq n$, are finitely generated, then $X$ is a generalized $n$-manifold.
\end{thm}
Another result related to the Modified Bing-Borsuk conjecture was obtained by Bryant \cite{br}:
\begin{thm}\cite{br}
Every $n$-dimensional homologically arc-homogeneous $ANR$-compactum is a generalized manifold.
\end{thm}

Here, a space $X$ is homologically arc-homogeneous \cite{qu} if for every path $\alpha:\mathbb I=[0,1]\to X$ the inclusion induced map
$$H_*(X\times\{0\},X\times\{0\}-(\alpha(0),0))\to H_*(X\times\mathbb I,(X\times\mathbb I)-\Gamma(\alpha))$$
is an isomorphism, where $\Gamma(\alpha)$ is the graph of $\alpha$.

More information about generalized manifolds can be found in Bryant \cite{br3}.

The last two theorems in this section show that $\mathbb Z$-homology manifolds have also common properties with Euclidean manifolds.
Recall that a space $X$ is a {\em Cantor $n$-manifold} \cite{hw}, \cite{ur} if  $X$ cannot be separated by a closed subset $F$ of dimension $\leq n-2$,
(i.e., $X\backslash F$ is disconnected).

\begin{thm}\cite{kr1}
Let $X$ be a locally compact, locally connected $\mathbb Z$-homology $n$-manifold with $\dim X=n>1$ at each point. Then $X$ is a local Cantor manifold, i.e. every open connected subset of $X$ is a Cantor $n$-manifold.
\end{thm}
This result was extended in \cite[Corollary 4.2]{tv}.
\begin{thm}\cite{tv} Let $X$ be a complete metric space which is a $\mathbb Z$-homology $n$-manifold. Then every open arcwise connected subset of $X$ is a Mazurkiewicz arc manifold with respect to the class of all spaces of dimension $\leq n-2$.
\end{thm}

Note that a space $X$ (not necessarily metrizable) is a {\em Mazurkiewicz arc manifold} with respect to the class of all spaces of dimension $\leq n-2$ \cite{tv}
if for ever two closed disjoint sets $A, B\subset X$, both having non-empty interiors in $X$, and every $F_\sigma$-set $F\subset X$ with $\dim F\leq n-2$, there is an arc $C$ in $X\backslash F$ joining $A$ and $B$.

\section{Homogeneous $ANR$-spaces}

We show in this section that finite-dimensional homogeneous $ANR$-spaces share many properties with Euclidean manifolds. In particular, the local cohomological and homological structure of homogeneous $n$-dimensional $ANR$-spaces is similar to the corresponding local structure of $\mathbb R^n$.
We also discuss another two problems of Bing-Borsuk \cite{bb} and their relation to the problem 
whether there exists a finite-dimensional non-degenerated homogeneous $AR$-compactum.

We say that a finite-dimensional space $X$ is {\em dimensionally full-valued} if $\dim X\times Y=\dim X+\dim Y$ for any compactum $Y$. 
It is known that all polyhedra and all one-dimensional compacta are dimensionally full-valued.
Pontryagin \cite{po} constructed in 1930 a family $\{\Pi_p:p{~}\mbox{is prime}\}$ of 2-dimensional homogeneous but not $ANR$-compacta such that $\dim (\Pi_p\times\Pi_q)=3$ for $p\neq q$. During the same time Borsuk raised the question whether $\dim X\times Y=\dim X+\dim Y$ for any $ANR$-compacta $X$ and $Y$. Kodama \cite{ko} provided a partial answer of Borsuk's question by proving that every 2-dimensional $ANR$-compactum is dimensionally full-valued. 
In 1988 Dranishnikov \cite{dra} gave a negative answer to Borsuk's question by constructing a family of $4$-dimensional metric $ANR$-compacta $M_p$, where $p$ is a prime number, such that $\dim (M_p\times M_q)=7$ for all $p\neq q$. The spaces $M_p$ are not homogeneous. After Dranishnikov constructed his examples, the question whether homogeneous $ANR$-compacta are dimensionally full-valued was raised. It goes back to \cite{br1} and was also discussed in \cite{clqr} and \cite{f}. 

This question was answered recently (for 3-dimensional homogeneous $ANR$-compacta it was known earlier \cite{vv1}).
\begin{thm}\cite{VV1}
Let $X$ be a finite-dimensional locally homogeneous $ANR$-space. Then the following holds:
\begin{itemize}
\item[(i)] $X$ is dimensionally full-valued;
\item[(ii)] If $X$ is homogeneous, then every $x\in X$ has a neighborhood $W_x$ such that $\bd\, \overline U$ is dimensionally full-valued for all $U\in\mathcal B_x$ with $\overline U\subset W_x$.
\end{itemize}
\end{thm}

According to \cite{bol}, a finite-dimensional compactum $X$ is dimensionally full-valued if and only if $\dim_GX=\dim X$ for any group $G$. It was shown in
\cite{vv1} that an $n$-dimensional $ANR$-compactum $X$ is dimensionally full-valued iff there exists a point $x\in X$ with 
$\check{H}_n(X,X\backslash x;\mathbb Z)$ is not trivial. 

Suppose $(K,A)$ is a pair of closed subsets of a space $X$ with $A\subset K$. Then we denote by $j^n_{K,A}:\check{H}^n(K;G)\to\check{H}^n(A;G)$
the inclusion induced cohomology homomorphism (recall that $dim_GX\leq n$ if and only if $j^n_{X,A}$ is surjective for every closed $A\subset X$).
We say that an element $\gamma\in\check{H}^n(A;G)$ is not extendable over $K$ if $\gamma$ is not contained in the image $j^n_{K,A}(\check{H}^n(A;G))$.
If $(K,A)$ is as above, $K$ is called an {\em $(n,G)$-cohomology membrane spanned on $A$ for an element  $\gamma\in\check{H}^n(A;G)$} if $\gamma$ is not extendable over $K$, but it is extendable over any proper closed subset $P$ of $K$ containing $A$.
The continuity of the \v{C}ech cohomology implies the following fact: If $A$ is a closed subset of a compact space $X$ and $\gamma\in\check{H}^n(A;G)$ is not extendable over $X$, then there is an $n$-cohomology membrane for $\gamma$ spanned on $A$. A space $X$ is said to be a {cohomological \em $(n,G)$-bubble} 
\cite{yo} if $\check{H}^n(X;G)\neq 0$ but $\check{H}^n(B;G)=0$ for every closed proper subset $B\subset X$.

The next theorem shows that the local cohomological structure of homogenous $n$-dimensional $ANR$-spaces is similar to the local structure of $\mathbb R^n$ (this was established earlier in \cite{vv1} for homogeneous $ANR$-compacta and countable principal domains $G$).
\begin{thm}\cite{VV1} Let $X$ be a connected homogeneous $ANR$-space with $\dim X=n\geq 2$ and $G$ be a countable group.
Then every point $x\in X$ has a basis $\mathcal B_x$ of open sets $U\subset X$ satisfying the following conditions:
\begin{itemize}
\item[(1)] $\rm{int}\overline U=U$ and the complement of $\bd\, U$ has exactly two components;
\item[(2)] $\check{H}^{n-1}(\bd\, U;G)\neq 0$, $\check{H}^{n-1}(\overline U;G)=0$ and $\overline{U}$ is an $(n-1,G)$-cohomology membrane spanned on $\bd\, U$ for any non-zero $\gamma\in\check{H}^{n-1}(\bd\, U;G)$;
\item[(3)] $\bd\, U$ is a cohomological $(n-1,G)$-bubble;
\end{itemize}
\end{thm}
A similar description of the local homology structure of homogeneous $ANR$-compacta is given in \cite{vv}.

We say that a space $X$ has an {\em $n$-dimensional $G$-obstruction} at a point $x\in X$ \cite{ku} if there is $W\in\mathcal B_x$  such that the homomorphism $j^n_{U,W}:H^{n}(X,X\backslash U;G)\to H^{n}(X,X\backslash W;G)$ is nontrivial
for every $U\in\mathcal B_x$  with $U\subset W$. Kuzminov \cite{ku} proved that every compactum $X$ with $\dim_GX=n$ contains a compact set $Y$ with $\dim_GY=n$ such that $X$ has an $n$-dimensional $G$-obstruction at any point of $Y$.

Theorems 3.1-3.2 provides more properties of homogeneous $n$-dimensional spaces which are typical for $n$-manifolds.
\begin{cor} Let $X$ be a locally homogeneous $ANR$-spaces with $\dim_X=n\geq 2$ and $G$ be a countable group. Then
\begin{itemize}
\item[(1)] $f(U)$ is open in $X$ provided $U\subset X$ is open and $f:U\to X$ is an injective map;
\item[(2)] $\dim A=n$, where $A\subset X$ is closed, if and only if $A$ has a non-empty interior in $X$;
\item[(3)] $\check{H}^{n}(P;G)=0$ for any proper compact set $P\subset X$;
\item[(4)] $X$ has an $n$-dimensional $G$-obstruction at every $x\in X$. Moreover, there is $W\in\mathcal B_x$ such that the   homomorphism $j_{U,V}^n$ is surjective for any $U,V\in\mathcal B_x$ with $\overline U\subset V\subset\overline V\subset W$.
\end{itemize}
\end{cor}
Properties (1) and (2) were also established by Lysko \cite{ly} and Seidel \cite{se}.
 
We say that $X$ is {\em cyclic in dimension $n$} if there is a group $G$ such that $\check{H}^n(X;G)\neq 0$. If a space is not cyclic in dimension $n$, it is called {\em acyclic in dimension $n$}. 
If $X$ is an $n$-dimensional $ANR$-compactum, the duality \cite{hw} between \v{C}ech homology and cohomology, and the universal coefficient formulas imply the following equivalence: $\check{H}_n(X;G)\neq 0$ for some group $G$ if and only if $X$ cyclic in dimension $n$.

We denote by $\mathcal H(n)$ the class of all homogeneous metric $ANR$-compacta of dimension $n$. 
\begin{qu}\cite{bb} Let $X\in\mathcal H(n)$. Is it true that:
\begin{itemize}
\item[(1)] $X$ is cyclic in dimension $n$$?$
\item[(2)] No closed subset of $X$, acyclic in dimension $n-1$, separates $X$$?$
\end{itemize}
\end{qu}
 
The next theorem shows that the two parts of  Question 3.4 have  positive or negative answers simultaneously.
\begin{thm}\cite{vv}
The following conditions are equivalent:
\begin{itemize}
\item[(1)] For all $n\geq 1$ and $X\in\mathcal H(n)$ there exists a group $G$ with $\check{H}^n(X;G)\neq 0$ $($resp., $\check{H}_n(X;G)\neq 0)$;
\item[(2)] If $X\in\mathcal H(n)$, $n\geq 1$, and $F\subset X$ is a closed set separating $X$, then there exists a group $G$ with $\check{H}^{n-1}(F;G)\neq 0$ $($resp., $\check{H}_{n-1}(F;G)\neq 0)$;
\item[(3 )] If $X\in\mathcal H(n)$, $n\geq 1$, and $F\subset X$ is a closed set separating $X$ with $\dim F\leq n-1$, then there exists a group $G$ such that $\check{H}^{n-1}(F;G)\neq 0$ $($resp., $\check{H}_{n-1}(F;G)\neq 0)$.
\end{itemize}
\end{thm}
Note that for any finite-dimensional homogeneous continuum $X$ (not necessarily $ANR$) we have the following result \cite{ktv}: If $\check{H}^n(X;G)\neq 0$ for some group $G$, then $\check{H}^{n-1}(F;G)\neq 0$ for any closed set $F\subset X$ separating $X$ with $\dim_GF\leq n-1$.  

On the other hand, the structure of cyclic homogeneous $ANR$ continua is described in \cite{vv3} (the notion of strong $V^n_G$-continua is given in  Section 4). 
\begin{thm}\cite{vv3}
Let $X$ be a homogeneous metric $ANR$-continuum such that $\dim_GX=n$ and $\check{H}^n(X;G)\neq 0$ for some group $G$. Then
\begin{itemize}
\item[(1)] $X$ is a cohomological $(n,G)$-bubble;
\item[(2)] $X$ is a strong $V^n_G$-continuum;
\item[(3)] $\check{H}^{n-1}(A;G)\neq 0$ for every closed set $A\subset X$ separating $X$.
\end{itemize}
\end{thm}
Items $(1)$ and $(3)$ were also established by Yokoi \cite{yo} for the case $G$ is a principal ideal domain.

Clearly, the cyclicity of finite-dimensional homogeneous $ANR$-compacta provides a negative answer to the next question.
\begin{qu}\cite{bb}, \cite{bo} Does there exists a non-degenerate finite-dimensional homogeneous $AR$-compactum$?$
\end{qu}

According to Fadell \cite{fa} there is no non-degenerate strongly homogeneous space. Here, a compactum $X$ is {\em strongly homogeneous} if $X$ is connected and for every $x_1\in X$ there is an open set $U$ containing $x_1$ and a continuous map $L_U:U\to C(X,X)$ with the following properties ($C(X,X)$ is the space of continuous maps from $X$ into $X$ with the compact-open topology): (i) For every $x\in U$, $L_U(x):X\to X$ is a homeomorphism such that $L_U(x)(x_1)=x$; (ii) $L_U(x_1)$ is the identity homeomorphism on $X$. This kind of homogeneity seems to be quite strong because Fadell's result implies that the Hilbert cube $Q$ is not strongly homogeneous. 

On the other hand it is interesting if we consider strongly locally homogeneous spaces in Question  3.7.
\begin{qu} Does there exists a non-degenerate finite-dimensional strongly locally homogeneous $AR$-compactum$?$
\end{qu}

Another questions in that direction was listed in \cite{west}.
\begin{qu}\cite{west} Is the Hilbert cube $Q$ the only homogeneous non-degenerate compact $AR$$?$
\end{qu}

\section{Separation of homogeneous spaces}

We already observed that the existence of finite-dimensional $AR$-compacta is equivalent to the question whether homogeneous $n$-dimensional $ANR$-compacta can be separated by closed sets $A$ acyclic in dimension $n-1$ with $\dim A\leq n-1$.   
In this section we discuss the question of separating homogeneous $n$-dimensional spaces by sets of a smaller dimension. Cantor manifolds defined in Section 2 (just before Theorem 2.6) were introduced by Urysohn \cite{ur} in 1925 as a generalization of Euclidean manifolds.
One of the first result concerning separation of homogeneous spaces was established by Krupski \cite{kr2}, \cite{kr3}.
\begin{thm}
Every region in a homogeneous $n$-dimensional space cannot be separated by a subset of dimension $\leq n-2$.
\end{thm}
The notion of Cantor manifolds was generalized in different ways.
Inspired by the classical result of Mazurkiewicz that any region in $\mathbb R^n$ cannot be cut by subsets of dimension $\leq n-2$
(a subset cuts if its compliment is not continuum-wise connected), Hadjiivanov-Todorov \cite{ht} introduced the class of {\em Mazurkiewicz manifolds}.
This notion was generalized in \cite{kktv} as follows:
A normal space (not necessarily metrizable) $X$  is a {\em Mazurkiewicz manifold with respect to  $\mathcal{C}$}, where
$\mathcal{C}$ is a class of spaces,
if for every two closed, disjoint subsets $X_0,X_1\subset X$, both
having non-empty interiors in $X$, and every $F_\sigma$-subset
$F\subset X$ with $F\in\mathcal{C}$, there exists a continuum $K$ in
$X\setminus F$ joining $X_0$ and $X_1$. Obviously, every Mazurkiewicz manifold with respect to the class at most $(n-2)$-dimensional spaces is
a Cantor $n$-manifold.

A new dimension $\mathcal D_{\mathcal K}$, unifying both the covering and the cohomological dimension, was introduced in \cite{kktv}. By $\mathcal D_{\mathcal K}^n$ we denote all spaces $X$ with $\mathcal D_{\mathcal K}(X)\leq n$. Concerning that dimension, we have the following result:
\begin{thm}\cite{kv} Let $X$ be a homogeneous locally connected space.
Then every region $U\subset X$ with $D_{\mathcal{K}}(U)=n$ is a Mazurkiewicz manifold with respect to the class $\mathcal D_{\mathcal K}^{n-2}$.
\end{thm}
Alexandroff \cite{ps} introduced another property which is possessed by compact closed $n$-manifolds, to so-called continua $V^n$. Here is the general notion of  Alexandroff manifold, see \cite{tv}:
A connected space $X$ is an {\em Alexandroff manifold with respect to a given class $\mathcal{C}$} of spaces if for every two disjoint closed subsets $X_0,X_1$ of $X$, both having non-empty interiors, there exists an open cover $\omega$ of $X$ such that no partition $P$ between $X_0$ and $X_1$ admits an $\omega$-map onto a space $Y\in\mathcal{C}$. The Alexandroff {\em continua $V^n$} are compact Alexandroff manifolds with respect to the class of all spaces $Y$ with $\dim Y\leq n-2$.
Recall that a partition between two disjoint sets $X_0,X_1$ in $X$ is a closed set $F\subset X$ such that $X\setminus F$ is the union of two open disjoint sets $U_0,U_1$ in $X$ with $X_0\subset U_0$ and $X_1\subset U_1$. An $\omega$-map $f:P\to Y$ is such a map that $f^{-1}(\gamma)$  refines $\omega$ for some open cover $\gamma$ of $Y$. 

A cohomological version of $V^n$-continua was considered in \cite{s}.
A compactum $X$ is a {\em $V^n_G$-continuum} \cite{s}, where $G$ is a given group, if for every open disjoint subsets $U_1,U_2$ of $X$ there is an open cover
$\omega$ of $X_0=X\setminus (U_1\cup U_2)$ such that any partition $P$ in $X$ between $U_1$ and $U_2$ does not admit an $\omega$-map $g$ onto a space $Y$ with $g^*:\check{H}^{n-1}(Y;G)\to \check{H}^{n-1}(P;G)$ being a trivial homomorphism. If, in addition, there is also an element
$\gamma\in\check{H}^{n-1}(X_0;G)$ such that for any partition $P$ between $U_1$ and $U_2$ and any $\omega$-map $g$ of $P$ into a space $Y$ we have $0\neq i_P^*(\gamma)\in g^*(\check{H}^{n-1}(Y;G))$, where $i_P$ is the embedding $P\hookrightarrow X_0$, $X$ is called a {\em strong $V^n_G$-continuum} \cite{vv3}.
Because $\check{H}^{n-1}(Y;G)=0$ for every space $Y$ with $\dim Y\leq n-2$, every $V^n_G$-continuum is a $V^n$-continuum in the sense of Alexandroff.

The following question, raised in \cite{tv}, is one of the remaining open problems concerning separation of homogeneous $ANR$-spaces.
\begin{qu} Let $X$ be a homogeneous $ANR$-continuum with $\dim X=n$ and $G$ be a group.
\begin{itemize}
\item[(1)] Is $X$ a $V^n$-continuum?
\item[(2)] Is $X$ a $V^n_G$-continuum?
\end{itemize}
\end{qu}
According to Theorem 3.6, Question 4.3 has a positive answer provide $\check{H}^{n}(X;G)\neq 0$. For strongly locally homogeneous spaces (not necessarily $ANR$s) the answer of Question $4.3(1)$ is also positive. 
\begin{thm}\cite{kktv1}
Every strongly locally homogeneous connected space with $\dim_GX=n$ is an Alexandroff manifold with respect to the class of spaces $Y$ with $\dim_GY\leq n-2$.
\end{thm}

Finally, let's mention another recent result extending Theorem 4.1, as well as, the result from \cite{vmv} that no region in 2-dimensional strongly locally homogeneous space cannot be separated by an arc (we say that a closed set $C\subset X$ is irreducibly separating $X$ if there are two disjoint open sets $G_1,G_2$ in $X$ such that $\overline G_1\cap\overline G_2=C$ and $X=\overline G_1\cup\overline G_2$). We say that a point $x\in X$ has a {\em special base} $\mathcal B_x$ if for any neighborhoods $U,V$ of $x$ in $X$ with $\overline U\subset V$ there is
$W\in\mathcal B_x$ such that $\rm{bd}W$ separates $\overline V\backslash U$ between $\rm{bd}\overline V$ and $\rm{bd}\overline U$. 
\begin{thm}\cite{v}
Let $\Gamma$ be a region in a finite-dimensional homogeneous space $X$ with $\dim_GX=n\geq 2$, where $G$ is a countable Abelian group. Then $\Gamma$ cannot be irreducibly separated by any closed set $C\subset X$ with the following property:
\begin{itemize}
\item[(i)] $\dim_G C\leq n-1$ and $H^{n-1}(C;G)=0$;
\item[(ii)] There is a point $b\in C\cap\Gamma$ having a special local base $\mathcal B_C^b$ in $C$ with $H^{n-2}(\rm{bd}_CU;G)=0$ 
for every $U\in\mathcal B_C^b$.
\end{itemize}
If $X$ is strongly locally homogeneous, the finite-dimensionality of $X$ can be omitted and condition $(ii)$ can be weakened to the following one:
\begin{itemize}
\item[(iii)]
There is $b\in C\cap\Gamma$ having an ordinary base $\mathcal B_C^b$ in $C$ with $H^{n-2}(\rm{bd}_CU)=0$, $U\in\mathcal B_C^b$.
\end{itemize}
\end{thm}  



\begin{thebibliography}{999}


\bibitem{ps} P.~Alexandroff, Die Kontinua $(V^p)$ - eine
Versch\"{a}rfung der Cantorschen Mannigfaltigkeiten, Monatshefte
fur Math. 61 (1957), 67--76.

\bibitem{ba}
J.~Bales, Representable and strongly locally homogeneous spaces and strongly n-homogeneous spaces, Houston J. Math 2, 3 (1976), 315--327.




\bibitem{bb}
R.~H.~Bing and K.~Borsuk, Some remarks concerning topological homogeneous spaces,  Ann. of Math. 81, 1 (1965), 100--111.

\bibitem{bo}
K.~Borsuk, Theory of retracts, Monografie Matematyczne 44, PWN, Warsaw, 1967.


\bibitem{bol}
V.~Boltyanskii, On dimensional full-valuedness of compacta, Dokl. Akad. Nauk SSSR, 67 (1949), 773--777 (in Russian).

\bibitem{bre}
G.~Bredon, Sheaf Theory, Sec. Ed., Graduate texts in Mathematics 170, Springer, 1997.

\bibitem{br}
J.~Bryant, Homologically arc-homogeneous $ENR$s, Geometry and Topology Monographs 9 (2006), 1--6.

\bibitem{br1}
J.~Bryant, Reflections on the Bing-Borsuk conjecture, Abstracts of the talks presented at the 19th Annual Workshop in Geometric Topology,
June 13-15, 2002, 2-3.

\bibitem{br3}
J.~Bryant, A survey of recent results on generalized manifolds, Topology Apll. 113 (2001), 13--22.

\bibitem{br2}
J.~Bryant, Homogeneous $ENR$'s, Topology Appl. 27 (1987), 301--306.


\bibitem{BF}
J.~Bryant and S.~Ferry, An alpha approximation theorem for homology manifolds, preprint.








\bibitem{can}
J.~Cannon, The recognition problem: what is a topological manifold, Bull. Amer. Math. Soc. 84 (1978), 832--866.

\bibitem{clqr}
M.~Cardenas, F.~Lasheras,A.~Quintero and D.~Repov\v{s}, On manifolds with nonhomogeneous factors, Centr. Eur. J. Math. 10 (2012), 857--862.





\bibitem{dra1}
A.~Dranishnikov, Cohomological dimension theory of compact metric spaces, Topology Atlas invited contribution, vol. 6, 2001, 7--73

\bibitem{dra}
A.~Dranishnikov. Homological dimension theory, Russian Math. Surveys 43, 4 (1988), 11--63.






\bibitem{fa}
E.Fadell, A note on the non-existence of strongly homogeous $AR$'s, Bull. Acad. Polon. Sci. Ser. Sci. Math. Astronom. Phys. 12 (1964), 531--534.

\bibitem{f}
V.~Fedorchuk, On homogeneous Pontryagin surfaces, Dokl. Akad. Nauk 404, 5 (2005), 601--603 (in Russian).

\bibitem{ja1}
W.~Jakobsche, Homogeneous cohomology manifolds which are inverse limits, Fund. Math. 137 (1991), 81--95.

\bibitem{ja}
W.~Jakobsche, The Bing-Borsuk conjecture is stronger than the Poincare conjecture, Fund. Math. 106 (1980), 127--134.

\bibitem{ht}
N.~Hadjiivanov and V.~Todorov, On non-Euclidean manifolds, C. R. Acad. Bulgare Sci. 33 (1980), 449--452 (in Russian).

\bibitem{hr}
D.~Halverson and D.~Repov\v{s},. The Bing-Borsuk conjecture and the Busemann conjecture, Math. Communications, 13 (2008), 163--184.


\bibitem{hw} W. Hurewicz and H. Wallman, Dimension theory.
Princeton University Press, Princeton, 1948.


\bibitem{kktv1}
A.~Karassev, P.~Krupski, V.~Todorov and V.~Valov, On generalized $V^n$-spaces, arXiv:2303.16373 [math.GN].

\bibitem{ktv}
A.~Karassev, V.~Todorov and V.~Valov, Alexandroff manifolds and homogeneous continua, Canad. Math. Bull. 57, 2 (2014), 335--343.

\bibitem{kktv}
A.~Karassev, P.~Krupski, V.~Todorov and V.~Valov, Generalized Cantor manifolds and homogeneity, Houston J. Math. 38, 2 (2012), 583--609.

\bibitem{ko}
Y.~Kodama, On homotopically stable points, Fund. Math. 44 (1957), 171–185.

\bibitem{kr1}
P.~Krupski, On the disjoint $(0,n)$-cells property for homogeneous ANRs, Colloq. Math. 66, 1 (1993), 77—84.

\bibitem{kr2}
P.~Krupski, Recent results on homogeneous curves and ANRs, Topology Proc. 16 (1991), 109—118.

\bibitem{kr3}
P.~Krupski, Homogeneity and Cantor manifolds, Proc. Amer. Math. Soc. 109 (1990), 1135--1142.

\bibitem{kv}
P.~Krupski and V.~Valov, Mazurkiewicz manifolds and homogeneity, Rocky J. Math. 41, 6 (2011), 1933--1938.

\bibitem{ku} V.~Kuz'minov, Homological dimension theory, Russian Math. Surveys 23 (1968), no. 1, 1--45.

\bibitem{ly}
J.~Lysko, Some theorems concerning finite dimensional homogeneous ANR-spaces, Bull. Acad. Polon. Sci. Sér. Sci. Math. Astronom. Phys.
24, 7 (1976), 491–496.

\bibitem{vmv}
J. van Mill and V. Valov, Homogeneous continua that are are not separated by arcs,
Acta Math. Hung. 157 (2019), 364--370.


\bibitem{po}
L.~Pontryagin, Sur une hypothese foundamentale de la dimension, C.R. Acad. Sci. 190 (1930), 1105--1107.

\bibitem{qu}
F.~Quinn. Problems on homology manifolds. \textit{Geometry and Topology Monographs} \textbf{9} (2006), 87--103.


\bibitem{se}
H.~Seidel, Locally homogeneous $ANR$-spaces, Arch. Math. 44 (1985), 79--81.


\bibitem{s}
S.~Stefanov, A cohomological analogue of $V^n$-continus and a theorem of Mazurkiewicz, Serdica Math. J. 12, 1 (1986), 88--94.


\bibitem{tv}
V.~Todorov and V.~Valov, Alexandroff type manifolds and homology manifolds, Houston J. Math. 40, 4 (2014), 1325--1346.

\bibitem{ur}
P.~Urysohn, Memoire sur les multiplicites cantoriennes, Fund. Math. 7 (1925), 30--137.

\bibitem{yo}
K.~Yokoi, Bubbly continua and homogeneity, Houston J. Math. 29, 2 (2003), 337--343.

\bibitem{v}
V.~Valov, Separation of homogeneous connected locally compact spaces, preprint.

\bibitem{VV1}
V.~Valov, Local structure of homogeneous $ANR$-spaces, arXiv:2303.10205v1 [math.GN].

\bibitem{VV}
V.~Valov, Homogeneous metric $ANR$-compacta, Serdica Math. Journal, 46 (2020), 1--18.

\bibitem{vv}
V.~Valov, Local homological properties and cyclicity of homogeneous $ANR$ compacta,
Proc. Amer. Math. Soc. 146 (2018), 2697--2705.



\bibitem{vv1}
V.~Valov, Local cohomological properties of homogeneous $ANR$ compacta, Fund. Math. 223 (2016), 257--270.

\bibitem{vv3}
V.~Valov, Homogeneous $ANR$-spaces and Alexandroff manifolds, Topology Appl. 173 (2014), 227--233.

\bibitem{west}
J.~West, Open problems in infinite-dimensional topology. in: Open Problems in Topology (J. van Mill and G.M. Reed, Eds.),
Elsevier Science Publishers B. V., North-Holland, 1990, 524--597.

\bibitem{w}
R.~Wilder, Topology of manifolds, Amer. Math. Soc. Coll. 32 (1949).
\end{thebibliography}
\end{document}